\documentclass{amsart}
\usepackage{graphicx} 
\usepackage{comment}

\author{Marni Mishna, Lucas Pannier, and Juan Pulido}
\address{MM and JP:~Department of Mathematics, Simon Fraser University, Burnaby, Canada\\ LP:~Université Paris-Saclay, UVSQ, CNRS, Laboratoire de mathématiques de Versailles, 78000, Versailles, France.}
\title{On Small-step quarter plane lattice walks with a non D-finite univariate generating function}

\usepackage{xcolor}
\definecolor{articleblue}{HTML}{16324F}
\definecolor{deepteal}{HTML}{0B6E6E}
\definecolor{pathred}{HTML}{B23A48}
\definecolor{softgold}{HTML}{B8872B}
\definecolor{linkblue}{HTML}{1A5FB4}
\definecolor{darkgray}{HTML}{444444}
\definecolor{asparagus}{rgb}{0.53,0.66,0.42}
\definecolor{purple1}{rgb}{0.45,0.32,0.62}
\definecolor{commentcolor}{RGB}{0,130,140} 
\usepackage{amsthm, amsfonts}
\usepackage{fullpage} 

\newtheorem{theorem}{Theorem}
\newtheorem{conjecture}{Conjecture}

\theoremstyle{definition}
\newtheorem{example}{Example}
\usepackage{amsmath,amssymb,amsfonts,amsthm,mathtools}
\usepackage{mathrsfs}
\usepackage{bm}
\usepackage{dsfont}
\usepackage{stmaryrd}
\usepackage{enumitem}

\usepackage{graphicx}
\usepackage{multirow}
\usepackage{float}
\usepackage{caption}
\usepackage{subcaption}
\usepackage{tikz}
\usepackage{tikz-cd}
\usetikzlibrary{
    arrows.meta,
    calc,
    decorations.pathmorphing,
    decorations.markings,
    intersections,
    matrix,
    patterns,
    positioning,
    shapes.geometric,
    backgrounds
}
\usepackage{pgfplots}
\pgfplotsset{compat=1.18}

\usetikzlibrary{arrows.meta}

\usepackage{longtable}
\usepackage{booktabs}
\usepackage{array}
\usepackage{tabularx}
\usepackage{makecell}
\usepackage{tikz}
\usetikzlibrary{arrows.meta}

\renewcommand{\arraystretch}{1.4}
\newcolumntype{M}[1]{>{\centering\arraybackslash}m{#1}}
\newcolumntype{L}[1]{>{\raggedright\arraybackslash}m{#1}}

\usepackage[
    colorlinks=true,
    linkcolor=deepteal,
    citecolor=pathred,
    urlcolor=linkblue,
    pdfauthor={Juan Pulido},
    pdftitle={Non-D-finiteness of a lattice-path generating function}
]{hyperref}
\usepackage[nameinlink,noabbrev]{cleveref}

\usepackage[
    backend=biber,
    style=numeric,
    maxbibnames=99,
    giveninits=true,
    doi=true,
    isbn=false,
    url=true
]{biblatex}

\newcommand{\StepModel}[1]{%
\begin{tikzpicture}[scale=0.4, baseline=-0.5ex]
    \draw[gray!45, rounded corners=1pt] (-1.35,-1.35) rectangle (1.35,1.35);

    \draw[gray!50, thin] (-1.15,0) -- (1.15,0);
    \draw[gray!50, thin] (0,-1.15) -- (0,1.15);

    \fill[black] (0,0) circle (1.6pt);

    \foreach \x/\y in {#1}{
        \draw[-{Latex[length=2mm]}, thick, articleblue]
            (0,0) -- (\x,\y);
    }
\end{tikzpicture}%
}
\newcommand{\DriftModel}[2][0.45]{%
\begingroup
    \def\driftx{0}%
    \def\drifty{0}%

    \foreach \x/\y in {#2}{%
        \pgfmathtruncatemacro{\newdriftx}{\driftx + \x}%
        \xdef\driftx{\newdriftx}%
        \pgfmathtruncatemacro{\newdrifty}{\drifty + \y}%
        \xdef\drifty{\newdrifty}%
    }%

    \pgfmathtruncatemacro{\DX}{\driftx}%
    \pgfmathtruncatemacro{\DY}{\drifty}%
    \pgfmathtruncatemacro{\DMAX}{max(abs(\DX),abs(\DY))}%

    \ifnum\DMAX=0
        \pgfmathsetmacro{\drawx}{0}%
        \pgfmathsetmacro{\drawy}{0}%
    \else
        \pgfmathsetmacro{\drawx}{\DX/\DMAX}%
        \pgfmathsetmacro{\drawy}{\DY/\DMAX}%
    \fi

    \begin{tikzpicture}[scale=0.5, baseline=-0.5ex]
        \draw[gray!45, rounded corners=1pt]
            (-1.35,-1.35) rectangle (1.35,1.35);

        \draw[gray!50, thin] (-1.15,0) -- (1.15,0);
        \draw[gray!50, thin] (0,-1.15) -- (0,1.15);

        \fill[black] (0,0) circle (1.6pt);

        \ifnum\DMAX=0
            \draw[pathred, thick] (0,0) circle (0.22);
            \fill[pathred] (0,0) circle (1.2pt);
        \else
            \draw[-{Latex[length=2mm]}, thick, pathred]
                (0,0) -- (\drawx,\drawy);
        \fi

        \node[
            font=\tiny,
            text=pathred,
            anchor=south east,
            inner sep=1pt
        ] at (1.28,-1.31) {$\scriptstyle(\DX,\DY)$};
    \end{tikzpicture}%
\endgroup
}

\newcommand{\DriftAndPolarDualCone}[4][1.25]{%
\begingroup
    \def\driftx{0}\def\drifty{0}%

    \foreach \x/\y in {#2}{%
        \pgfmathtruncatemacro{\newdriftx}{\driftx + \x}%
        \xdef\driftx{\newdriftx}%
        \pgfmathtruncatemacro{\newdrifty}{\drifty + \y}%
        \xdef\drifty{\newdrifty}%
    }%

    \pgfmathtruncatemacro{\DX}{\driftx}%
    \pgfmathtruncatemacro{\DY}{\drifty}%
    \pgfmathtruncatemacro{\DMAX}{max(abs(\DX),abs(\DY))}%

    \ifnum\DMAX=0
        \pgfmathsetmacro{\drawx}{0}%
        \pgfmathsetmacro{\drawy}{0}%
    \else
        \pgfmathsetmacro{\drawx}{0.92*\DX/\DMAX}%
        \pgfmathsetmacro{\drawy}{0.92*\DY/\DMAX}%
    \fi

    \begin{tikzpicture}[scale=0.5, baseline=-0.5ex]

        \path (-1.45,-1.45) rectangle (1.45,1.45);

        \draw[gray!35, thin] (-1.18,0) -- (1.18,0);
        \draw[gray!35, thin] (0,-1.18) -- (0,1.18);

        \begin{scope}
            \clip[rounded corners=1.5pt]
                (-1.32,-1.32) rectangle (1.32,1.32);

            \fill[asparagus!18]
                (0,0) -- ($#1*#3$) -- ($#1*#4$) -- cycle;

        \end{scope}

    \ifnum\DMAX=0
    \draw[pathred!65, line width=0.55pt] (0,0) circle (0.2);
    \fill[pathred!65] (0,0) circle (1.0pt);
\else
    \draw[
        -{Latex[length=1.9mm]},
        line width=0.6pt,
        pathred!70
    ]
        (0,0) -- ({0.95*\drawx},{0.95*\drawy});
\fi

        \fill[black] (0,0) circle (1.35pt);

        \draw[gray!40, rounded corners=1.5pt]
            (-1.32,-1.32) rectangle (1.32,1.32);

    \end{tikzpicture}%
\endgroup
}

\newcommand{\DriftAndRegion}[6][1.25]{%
\begingroup
    \def\driftx{0}\def\drifty{0}%

    \foreach \x/\y in {#2}{%
        \pgfmathtruncatemacro{\newdriftx}{\driftx + \x}%
        \xdef\driftx{\newdriftx}%
        \pgfmathtruncatemacro{\newdrifty}{\drifty + \y}%
        \xdef\drifty{\newdrifty}%
    }%

    \pgfmathtruncatemacro{\DX}{\driftx}%
    \pgfmathtruncatemacro{\DY}{\drifty}%
    \pgfmathtruncatemacro{\DMAX}{max(abs(\DX),abs(\DY))}%

    \ifnum\DMAX=0
        \pgfmathsetmacro{\drawx}{0}%
        \pgfmathsetmacro{\drawy}{0}%
    \else
        \pgfmathsetmacro{\drawx}{0.92*\DX/\DMAX}%
        \pgfmathsetmacro{\drawy}{0.92*\DY/\DMAX}%
    \fi

    \begin{tikzpicture}[scale=0.4, baseline=-0.5ex]

        \path (-1.45,-1.45) rectangle (1.45,1.45);

        \draw[gray!35, thin] (-1.18,0) -- (1.18,0);
        \draw[gray!35, thin] (0,-1.18) -- (0,1.18);

        \begin{scope}
            \clip[rounded corners=1.5pt]
                (-1.32,-1.32) rectangle (1.32,1.32);

          \draw[
    line width=0.85pt,
    asparagus!85!black
]
    (0,0) -- ($#1*#3$);

\draw[
    line width=0.85pt,
    asparagus!85!black
]
    (0,0) -- ($#1*#4$);

\draw[
    line width=0.85pt,
    asparagus!85!black
]
    (0,0) -- ($#1*#5$);

\draw[
    line width=0.85pt,
    asparagus!85!black
]
    (0,0) -- ($#1*#6$);
        \end{scope}

    \ifnum\DMAX=0
    \draw[pathred!65, line width=0.55pt] (0,0) circle (0.2);
    \fill[pathred!65] (0,0) circle (1.0pt);
\else
    \draw[
        -{Latex[length=1.9mm]},
        line width=0.6pt,
        pathred!70
    ]
        (0,0) -- ({0.95*\drawx},{0.95*\drawy});
\fi

        \fill[black] (0,0) circle (1.35pt);

        \draw[gray!40, rounded corners=1.5pt]
            (-1.32,-1.32) rectangle (1.32,1.32);

    \end{tikzpicture}%
\endgroup
}

\newcommand{\DriftAndRegiontemp}[6][1.25]{%
\begingroup
    \def\driftx{0}\def\drifty{0}%

    \foreach \x/\y in {#2}{%
        \pgfmathtruncatemacro{\newdriftx}{\driftx + \x}%
        \xdef\driftx{\newdriftx}%
        \pgfmathtruncatemacro{\newdrifty}{\drifty + \y}%
        \xdef\drifty{\newdrifty}%
    }%

    \pgfmathtruncatemacro{\DX}{\driftx}%
    \pgfmathtruncatemacro{\DY}{\drifty}%
    \pgfmathtruncatemacro{\DMAX}{max(abs(\DX),abs(\DY))}%

    \ifnum\DMAX=0
        \pgfmathsetmacro{\drawx}{0}%
        \pgfmathsetmacro{\drawy}{0}%
    \else
        \pgfmathsetmacro{\drawx}{0.92*\DX/\DMAX}%
        \pgfmathsetmacro{\drawy}{0.92*\DY/\DMAX}%
    \fi

    \begin{tikzpicture}[scale=0.5, baseline=-0.5ex]

        \path (-1.45,-1.45) rectangle (1.45,1.45);

        \draw[gray!35, thin] (-1.18,0) -- (1.18,0);
        \draw[gray!35, thin] (0,-1.18) -- (0,1.18);

        \begin{scope}
            \clip[rounded corners=1.5pt]
                (-1.32,-1.32) rectangle (1.32,1.32);

         \end{scope}

    \ifnum\DMAX=0
    \draw[pathred!65, line width=0.55pt] (0,0) circle (0.2);
    \fill[pathred!65] (0,0) circle (1.0pt);
\else
    \draw[
        -{Latex[length=1.9mm]},
        line width=0.6pt,
        pathred!70
    ]
        (0,0) -- ({0.95*\drawx},{0.95*\drawy});
\fi

        \fill[black] (0,0) circle (1.35pt);

        \draw[gray!40, rounded corners=1.5pt]
            (-1.32,-1.32) rectangle (1.32,1.32);

    \end{tikzpicture}%
\endgroup
}

\newcommand{\Qxy}{Q(x,y;t)}
\newcommand{\Qx}{Q(1,0;t)}

\newcommand{\Qzero}{Q(0,0;t)}
\newcommand{\Qtot}{Q(1,1;t)}

\newcommand{\mjm}[1]{{\large \sl \color{orange} #1}}

\begin{document}
\begin{abstract}
We report on the status of the conjecture of Bousquet-M\'elou and Mishna that the univariate counting generating function of a small-step quarter-plane lattice model is D-finite if and only if a particular group associated to the walk is finite. While the finite-group case is fully resolved, the infinite-group case remains incomplete. 
First we survey the known cases.  By invoking asymptotics results from probability theory and a case analysis completed by Bostan, Raschel and Salvy (2014), we establish that the generating functions of 14 infinite-group models are not D-finite, and that for an additional 18 models the generating function for walks that end on an axis is not D-finite. We conclude with a discussion of limitations of asymptotic arguments to resolve the remaining cases.
 \end{abstract}
\maketitle


\tableofcontents

\section{Introduction}
\label{sec:intro}

A small-step quarter-plane lattice model is a combinatorial class defined by its step set~$S$, a finite subset of $\{-1, 0, 1\}^2\setminus \{(0,0)\}$. Such a class is the set of all walks on the lattice ${\mathbb{N}}^2$ starting at the origin, taking steps from $S$ and never stepping out of the first quadrant~$\mathbb{N}^2$. More precisely, we call a sequence $\omega= \omega_1\omega_2\dots \omega_n$ with $\omega_i\in S$ a walk, and denote by~$|\omega|$ the number of steps, $n$. The endpoint of the walk is the vector sum of the steps. It is convenient to include one empty walk of length 0 in the class. The \emph{counting generating function} of this class is $Q_S(t):=\sum_{\omega} t^{|\omega|}$ where the sum ranges over all valid walks $w$. 

We address the following:   
\begin{quote}\sl
Which features of a quarter plane lattice path model determine the type of equations satisfied by the associated generating function $Q_S(t)=\sum_{\omega} t^{|\omega|}$?
\end{quote}

This question has generated remarkable interest over the past quarter century for many reasons, not only for the variety of successful approaches, but also for what these approaches have revealed about connections between points of view on series solutions of functional equations and combinatorics. We are motivated by the fact not only that combinatorial properties of small-step quarter-plane walks are correlated with properties of series, but that many of these connections can be established through general arguments, rather than on a case-by-case basis. An important example is the systematic study of the series $Q_S(t)$ undertaken by Bousquet-Mélou and Mishna~\cite{bousquet-melou_walks_2010}. Starting from the $2^{8}=256$ different models, they first settle the trivial models and those equivalent to walks confined to a half-line or half-plane, showing that in these cases $Q_S(t)$ is algebraic. After identifying models related by reflection across the diagonal $x=y$, 79 non-trivial, non-isomorphic models remain to be considered. They associate a \emph{group} to each of these models and show that 23 have finite group and 56 have infinite group, see Section~\ref{sec:background} for the definition. This group originates in the work of Fayolle, Iasnogorodski and Malyshev~\cite{fayolle_random_1999}. The finite/infinite group dichotomy leads to the following conjecture, whose current status we examine here.

\begin{conjecture}[Bousquet-M\'elou and Mishna, 2010~\cite{bousquet-melou_walks_2010}; Mishna 2009~\cite{mishna_classifying_2009}]
\label{conj:q11}
For any one of the 79 nontrivial small-step quarter-plane lattice path models, the counting generating function satisfies a linear differential equation with polynomial coefficients, that is, is D-finite, if and only if a particular group associated to the model is finite.
\end{conjecture}

 The \emph{endpoint generating function} is defined as the following series\footnote{The overloading of the name $Q$ is not ideal, but we use the multi-variable version sparingly, so hopefully the double duty of $Q$ is not too confusing.} in ${\mathbb{N}}[x,y][[t]]$: 
 \[
  Q_S(x,y; t) := \sum_{n=0}^\infty \sum_{(i,j)\in{\mathbb{N}}^2} a_S(i,j;n) x^iy^j t^n\]
where $a_S(i,j;n)$ is the number of walks with~$n$ steps that end at the point~$(i,j)$. We write $c_S(n):=\sum_{(i,j) \in \mathbb{N}^{2}} a_S(i,j;n).$
The univariate counting generating function is an evaluation of this series: $Q_S(t)=\Qtot=\sum_{n\geq0}c(n)t^{n}$.    Other evaluations are meaningful: $Q_S(0,0;t)$, $Q_S(1,0;t)$ and $Q_S(0,1;t)$ are the generating functions for the subset of walks that end at the origin, on the~$x$-axis, and on the~$y$-axis, respectively. Henceforth, when the step set $S$ is clear, we omit the subscript $S$ from the notation. 

One side of the conjecture is now well established: if the group of a model is finite, the counting generating function is D-finite. This was completed by Bousquet-M\'elou and Mishna, with an important case proved by Bostan and Kauers~\cite{bostan_complete_2010}. Furthermore, these works show that the multivariate version is true: if the group is finite, $\Qxy$ is D-finite with respect to $x$, $y$ and $t$. The more difficult result is the converse, which is not proved in the univariate case, but is proved in the multivariate case:
\begin{theorem}[Hardouin, 2026+, Theorem 3.3~\cite{hardouin_d-finiteness_2025}; Dreyfus, Elvey Price, Raschel 2026+~\cite{dreyfus_enumeration_2024}]
\label{thm:endpointnotdelfinite}  
If the group associated to a small-step quarter plane walk model is infinite, then the generating series~$\Qxy$ does not satisfy a linear differential equation in any of $x$, $y$ or $t$ with coefficients in $\mathbb{Q}(x,y,t)$. 
\end{theorem}

The two recent complete proofs build on results developed over the past decade, but use different approaches. Although this is a strong result, it does not, \emph{a priori}, prove the conjecture, since the specialization $x=y=1$ could simplify the series and change its nature. 
This does not happen for the specialization $x=y=0$. 
\begin{theorem}[Bostan, Raschel and Salvy, 2014~\cite{bostan_non-d-finite_2014}]
For each of the 56 small-step models with an infinite group, the associated series $\Qzero$ is either exactly 1 or is not D-finite. 
\end{theorem}
Although this result does not directly give information about $Q(1,1;t)$, an intermediate result from its proof is essential to our arguments; we recall it as Theorem~\ref{thm:irrational} below.

The classification of series---and, more generally, functions---according to the types of equations they satisfy has almost a hundred and fifty years history within mathematics. Series arising from combinatorics bring additional intuition, while the analysis of generating functions can be useful for combinatorial classifications. The counting generating function of a quarter-plane model belongs to $\mathbb{N}[[t]]$, and we consider whether it is rational (elements of $\mathbb{C}(t)$), algebraic (a solution to a non-trivial polynomial equation over $\mathbb{C}(t)$), D-finite (a solution to a linear differential equation with coefficients in $\mathbb{C}(t)$), differentially algebraic (the derivatives satisfy a polynomial equation over $\mathbb{C}(t)$). If a series is not differentially algebraic, it is differentially transcendental\footnote{also called hyper-transcendental}. Each of these classes fit into a containment hierarchy 
$$\text{rational}\subseteq \text{algebraic}\subseteq \text{D-finite}\subseteq \text{differentially algebraic}.$$

\subsection{Plan of the paper}
We begin in the next section recalling definitions and preliminary results. It is followed by a survey of the known results. In Section~\ref{sec:singular} we cite the results of Mishna and Rechnitzer~\cite{mishna_two_2009} and of Melczer and Mishna~\cite{melczer_singularity_2014} establishing that the generating function of each of the five so-called `singular` or genus 0 models is not D-finite. These models admit no nonempty excursions and the corresponding endpoint generating functions each satisfy $Q(0,0;t)=1$. This is followed by a discussion on the differentially transcendental models. 

In Section~\ref{sec:reluctant} we define the drift regime of a model, following the lead of Bostan \emph{et al.}~\cite{bostan_non-d-finite_2014} and associating to each model a cone. We combine intermediary results of~\cite{bostan_non-d-finite_2014} with probabilistic results on exit time of random walks of Denisov and Wachtel~\cite{denisov_random_2015} and Duraj~\cite{duraj_random_2014} to prove that the walks with infinite group with drift either zero or sufficiently negative all have non-D-finite generating functions. When Bousquet-M\'elou et al. remarked in~\cite{bernardi_counting_2021} that 15 non-singular models had non-D-finite counting series, they are likely alluding to the argument we develop here\footnote{We list two fewer models as two cases are periodic.}.

 Our first main result is the following, proved in Theorems~\ref{thm:zerodriftnondfinite} and~\ref{thm:interiordriftnondfinite}. This result is alluded to in~\cite{bernardi_counting_2021}, and here we list the models and give an explicit proof. 
\begin{theorem}\label{thm:main1}
The counting series for each of the following models is not D-finite:
\begin{center}
\StepModel{-1/-1,-1/0,0/1,1/-1,1/1}
%
\StepModel{-1/-1,-1/1,0/1,1/0,1/-1} 
%
\StepModel{-1/0,-1/1,0/-1,1/1,1/-1} 

\StepModel{-1/0,0/1,0/-1,1/-1} 
%
\StepModel{-1/0,0/1,-1/-1,1/-1} 
%
\StepModel{-1/-1,-1/0,0/-1,1/1} 
%
\StepModel{-1/0,0/1,0/-1,1/0,-1/-1}  
%
\StepModel{-1/-1,-1/0,0/1,0/-1,1/-1} 
%
\StepModel{-1/0,-1/1,0/1,0/-1,1/-1} 
%
\StepModel{-1/-1,-1/0,-1/1,0/1,1/-1} 
%
\StepModel{-1/-1,-1/0,-1/1,0/1,1/0,1/-1} 
\StepModel{-1/-1,-1/0,-1/1,0/1,0/-1,1/-1} 
\StepModel{-1/-1,-1/0,-1/1,0/-1,1/1,1/-1} 
%
\StepModel{-1/-1,-1/0,-1/1,0/1,0/-1,1/-1,1/0} 

\end{center}
\end{theorem}

In Section~\ref{sec:bdrynondfinite} we present our second result as Theorem~\ref{thm:boundary-series-nondfinite}. To our knowledge, this result is new.
\begin{theorem}\label{thm:main2}
For each of the following models, the counting series for walks that end on the $x$-axis, $Q(1,0;t)$, is not D-finite:
\begin{center}
\StepModel{-1/-1,1/0,-1/0,0/1} 
\StepModel{-1/-1,1/0,0/1,-1/1}
\StepModel{-1/1,-1/0,0/1,1/-1}  
\StepModel{-1/0,-1/1,0/-1,1/1} 
\StepModel{0/-1,1/0,-1/0,0/1,-1/1} 
\StepModel{0/1,1/1,-1/0,-1/-1,0/-1}
\StepModel{-1/1,-1/0,0/1,0/-1,1/1} 
\StepModel{-1/-1,1/0,-1/0,0/1,-1/1} 
\StepModel{-1/-1,-1/0,-1/1,0/-1,1/1} 
\StepModel{-1/-1,1/0,-1/0,1/1,0/1,-1/1}
\StepModel{-1/-1,-1/0,-1/1,0/1,0/-1,1/1}
\StepModel{-1/-1,0/-1,1/0,-1/0,1/1,-1/1} 
\StepModel{-1/-1,-1/0,-1/1,0/1,1/1,1/-1} 
\StepModel{-1/-1,-1/0,-1/1,0/1,0/-1,1/1,1/0} 
\StepModel{-1/0,0/1,0/-1,1/1}
\StepModel{-1/1,-1/-1,0/1,1/1,1/0} 
\StepModel{-1/0,0/1,-1/-1,1/1,1/0} 
\StepModel{-1/-1,0/-1,1/0,1/1,0/1,-1/1} 
\end{center}
\end{theorem} 

The models in both theorems are characterized by where their drift lies relative to the associated cone. We conclude with a table summarizing the conjectured singular exponents for each drift regime (defined below). This provides a companion to the table of asymptotics of Fayolle and Raschel~\cite{fayolle_some_2012}.

The numbers of the models refer to the numbering used in ~\cite{bousquet-melou_walks_2010} and ~\cite{bostan_non-d-finite_2014}. The OEIS numbers refer to entries in the Online Encyclopedia of Integer Sequences. 

%
\begin{table}
  \begin{tabular}{rlll} \toprule
    $\#$ & Type & Nature & (possible) Reasoning \\ \midrule \hline
    5 & singular models  & non D-finite &infinite singularities \\
    14 & polar or zero drift regime & non D-finite & asymptotics  \\
18 & & likely non D-finite & asymptotics of return to axis\\ 
   \midrule
    9\footnote{six overlap with above cases} & decouplable & differentially algebraic &   \\
     \bottomrule
    \end{tabular}

\caption{Summary of the infinite-group models studied}
\label{tab:Summary}
\end{table}

\section{Background and Notation}
\label{sec:background}

Let~$S$ be a fixed step set. The \emph{inventory} $P$ and associated \emph{kernel} $K$ are defined, respectively, by
\[P(x,y):= \sum_{(i,j)\in S}x^iy^j\qquad \text{and }\qquad K(x,y,t):=xy(1-tP(x,y)).\]
 
For small-step models starting at the origin and confined to the quarter plane, the endpoint generating function \(Q(x,y;t)\) satisfies a well-known functional equation, which follows from a simple combinatorial recurrence:  a walk is either the trivial walk or a shorter walk followed by a step. This translates into {\small
\begin{equation}\label{eq:functional-equation}
K(x,y,t)Q(x,y;t)
=
xy
-K(x,0,t)Q(x,0;t)
-K(0,y,t)Q(0,y;t)
+K(0,0,t)Q(0,0;t).
\end{equation}    }
We make use of the specialization of this equation at $x=y=1$:
\[(1-|\mathcal S|t)Q(1,1;t)=1-K(1,0,t)Q(1,0;t)-K(0,1,t)Q(0,1;t)+K(0,0,t)Q(0,0;t).\]

We now recall the definition of the \emph{group associated with $S$}. 
We write the inventory as
\[ P(x,y)=A_{-1}(x)y^{-1}+A_0(x)+A_{1}(x)y=B_{-1}(y)x^{-1}+B_0(y)+B_1(y)x.\]

For the non-trivial models considered here, the four polynomials
$A_{-1}(x),A_1(x),B_{-1}(y),B_1(y)$ are nonzero.
The \emph{group of the walk} associated to $S$, denoted  $G(S)$, introduced 
in~\cite{fayolle_random_1999} and imported into combinatorics 
in~\cite{bousquet-melou_walks_2010},  is the group of birational transformations generated by the two involutions 
\[
\Phi(x,y)
\coloneqq
\left(x^{-1}\frac{B_{-1}(y)}{B_1(y)},\,y\right),
\qquad
\Psi(x,y)
\coloneqq
\left(x,\,y^{-1}\frac{A_{-1}(x)}{A_1(x)}\right).
\]
These two transformations arise naturally as the involutions that swap the two 
roots of the kernel $K(x,y,t)$ viewed as a polynomial in $x$ (respectively $y$),
and each has order two. Their key property is that they preserve the step 
inventory:
\[
P(\Phi(x,y))=P(x,y),
\qquad
P(\Psi(x,y))=P(x,y).
\]
Hardouin and Singer showed that the order of the group is bounded by 12~\cite{hardouin_differentially_2021}, hence by iterating the composition $\Phi(\Psi(x,y))$ sufficiently many times, one can determine whether or not the group is finite. 

\subsection{D-finite functions}
A univariate function or series $F(t)$ is \emph{D-finite} if it satisfies a linear differential equation with polynomial coefficients. That is, there exist polynomials $p_0(t),\ldots,p_k(t)$, not all zero, such that \[p_0(t)F(t)+p_1(t)F^{(1)}(t)+\ldots+p_k(t)F^{(k)}(t)=0.\] A univariate D-finite function has only a finite number of singularities in the complex plane. A multi-variable series $F(x,y;t)$ in $\mathbb{Q}[x,y][[t]]$ is said to be D-finite if it is D-finite with respect to each of the variables $x,y$ and $t$, with coefficients in $\mathbb{Q}[x,y,t]$.

When the coefficient sequence of $Q(t)$ is asymptotically equivalent to $\kappa\rho^{-n}n^{\alpha}$, we say that $\alpha$ is the \emph{singular exponent}. If $Q(t)$ is D-finite and has this asymptotic form, then $\alpha$ is rational. 
The rationality of this exponent in D-finite functions follows from Theorem 19.1 of~\cite[p.111]{wasow_asymptotic_1965}, but here we use the
 formulation of Bostan et al.\cite[Theorem 3]{bostan_non-d-finite_2014} recalled below. For additional details, we refer the reader to their discussion. 
\begin{theorem}[Bostan et al. 2014, Theorem 3~\cite{bostan_non-d-finite_2014}]
\label{thm:alphaisirrational}
Let $(f_n)_{n\geq 0}$ be an integer valued sequence whose $n${th} term $f_n$ behaves asymptotically like $\kappa\rho^{-n}n^\alpha$ for some real constant $\kappa>0$. If the growth constant $\rho$ is transcendental or if the singular exponent $\alpha$ is irrational, then the generating function $F(t)=\sum_{n\geq 0} f_n t^n$ is not D-finite.
\end{theorem}

\subsection{The singular models}
\label{sec:singular}
There are five singular models characterized by the property that the step set is contained in a half plane. These models are shown in Figure~\ref{fig:singular}. The non-D-finiteness of $Q_S(t)$ was initially established in a case by case analysis~\cite{mishna_two_2009, melczer_singularity_2014} demonstrating that the generating function had an infinite number of singularities, a property inconsistent with D-finiteness. 

\begin{figure}
\StepModel{-1/1,0/1,1/-1}\hspace{1mm}\StepModel{-1/1,1/1,1/-1}\hspace{1mm}\StepModel{-1/1,0/1,1/0,1/-1}\hspace{1mm}\StepModel{-1/1,0/1,1/1,1/-1}\hspace{1mm}\StepModel{-1/1,0/1,1/1,1/0,1/-1}
\caption{The five singular models. From left to right, these are
    models \(1,2,13,15,\) and \(22\).}
\label{fig:singular}
\end{figure}

\begin{theorem}[Mishna and Rechnitzer, 2009 \cite{mishna_two_2009}; Melczer and Mishna, 2014\cite{melczer_singularity_2014}]
 All five singular small-step models with an infinite group have a non-D-finite counting series.
\end{theorem} 

Subsequently, it was established that the endpoint generating series is in fact differentially transcendental by framing $Q(x,0;t)$ and $Q(0,y;t)$ as solutions to a $q$-equation. Singular models have the property that solutions to the kernel polynomial $K(x,y,t)=0$ form a curve of genus 0, and hence has a rational parametrization. This fact gives rise to a relatively straightforward yet elegant proof which connects to Bernoulli numbers among other things. The connection to Bernoulli numbers is developed in~\cite{bostan_differential_2024}. However, they note ``that the nature of \(Q(1,1;t)\) still remains out of reach with these techniques".

\subsection{Differentially algebraic models}
\label{sec:dalg}
Among the infinite-group models, some have endpoint generating functions $Q(x,y;t)$ that are \emph{differentially algebraic},\footnote{A multivariate series $F(x,y;t)$ is \emph{differentially algebraic in $t$} if there exist an integer $r\geq 0$ and a nonzero polynomial 
$\Phi \in \mathbb{Q}(x,y,t)[z_0,z_1,\ldots,z_r]$ such that
\[
\Phi\bigl(F, \partial_t F, \ldots, \partial_t^r F\bigr)=0,
\]
and analogously for $x$ and $y$.} this is established by the existence of a
\emph{decoupling function}. For walks starting at the origin, this means that
there exist rational functions
\[
F(x) \in \mathbb{Q}(x,t)
\qquad \text{and} \qquad
G(y) \in \mathbb{Q}(y,t)
\]
such that, on the kernel curve \(K(x,y;t)=0\),
\[
xy = F(x)+G(y).
\]
This identity separates the mixed term \(xy\) into a part depending only on
\(x\) and a part depending only on \(y\), and is the key mechanism that allows
the boundary terms in the functional equation to be treated separately.
Bernardi, Bousquet-Mélou, and Raschel~\cite{bernardi_counting_2021} proved that, among the infinite-group
models, exactly nine are decoupled for walks starting at the origin.

\begin{theorem}[Bernardi, Bousquet-M\'elou and Raschel, 2021 \cite{bernardi_counting_2021}]
For each of the nine infinite-group models shown in
Figure~\ref{fig:dalg-models}, the endpoint generating function \(Q(x,y;t)\)
is differentially algebraic in \(x\), \(y\), and \(t\). Consequently, the
counting series \(Q(1,1;t)\) is differentially algebraic over $\mathbb{C}(t)$.
\label{thm:DAmodels}
\end{theorem}
\begin{figure}
    \centering

    \begin{tabular}{@{}ccccc@{}}
        \StepModel{-1/-1,-1/0,0/1,1/0} &
        \StepModel{-1/-1,-1/1,0/1,1/0} &
        \StepModel{-1/0,0/-1,0/1,1/1} &
        \StepModel{-1/0,0/1,1/0,1/-1} &
        \StepModel{-1/-1,-1/0,0/1,1/0,1/1}
        \\
        {\small Model 3} &
        {\small Model 4} &
        {\small Model 6} &
        {\small Model 8} &
        {\small Model 25}
    \end{tabular}

    \medskip

    \begin{tabular}{@{}cccc@{}}
        \StepModel{-1/-1,-1/0,0/1,0/-1,1/1} &
        \StepModel{-1/-1,-1/0,-1/1,0/1,1/0} &
        \StepModel{-1/0,0/1,0/-1,1/0,1/1} &
        \StepModel{-1/1,0/1,0/-1,1/0,1/-1}
        \\
        {\small Model 26} &
        {\small Model 37} &
        {\small Model 38} &
        {\small Model 39}
    \end{tabular}

    \caption{The nine infinite-group models that are decoupled for
    walks starting at the origin.}
    \label{fig:dalg-models}
\end{figure}
\noindent Interestingly, the existence of a decoupling function is not determined by the step set: it may also depend on the starting point. Bernardi, Bousquet-Mélou, and Raschel \cite{bernardi_counting_2021} show that in addition to the nine infinite-group models decoupled for walks starting at the origin, three further infinite-group models admit decoupling functions at other starting points; these models are shown in Figure~\ref{fig:dalg-other-starting-points}. The corresponding generating functions are expected to be differentially algebraic, although this is not proved there. 

\begin{figure}[h!]
    \begin{tabular}{@{}ccc@{}}
    \StepModel{-1/-1,-1/1,1/0,1/1} &
    \StepModel{-1/0,0/1,1/1,1/-1} &
    \StepModel{-1/-1,0/1,0/-1,1/-1,1/1}
    \\
    {\small Model 7} &
    {\small Model 14} &
    {\small Model 32}
\end{tabular}

    \caption{The three additional infinite-group models that admit
    decoupling functions for certain points other than the origin.}
    \label{fig:dalg-other-starting-points}
\end{figure}
Some models in Figure~\ref{fig:dalg-models} and \ref{fig:dalg-other-starting-points} are displayed up to reflection across
the diagonal \(x=y\). This preserves \(Q(1,1;t)\) and \(Q(0,0;t)\), while
interchanging the boundary series \(Q(1,0;t)\) and \(Q(0,1;t)\).

Returning to walks starting at the origin, Theorem~\ref{thm:DAmodels} does
not settle the \(D\)-finiteness of the total counting series \(Q(1,1;t)\):
every \(D\)-finite series is differentially algebraic, but the
converse does not hold in general. These nine models are therefore especially interesting: their infinite group suggest that \(Q(1,1;t)\) should be non-\(D\)-finite, while their decoupling functions imply differential algebraicity. The additional structure of these models even raised the possibility that $Q(1,1;t)$ might be algebraic.
\section{The drift regime  and the probability framework}
\label{sec:reluctant}

Bostan, Raschel and Salvy~\cite{bostan_non-d-finite_2014} use asymptotic results for random walks in cones to obtain an explicit expression for the singular exponent $\alpha$ for $Q(0,0;t)$. For each of the 51 non-singular infinite group models, they then show by case analysis that $\alpha$ is irrational and hence  $Q(0,0;t)$ is not D-finite.  

They use a probabilistic argument appealing to results of Denisov and Wachtel~\cite{denisov_random_2015}. The main idea in our work is very similar: we apply related results of Denisov and Wachtel~\cite{denisov_random_2015}, and also results of Duraj~\cite{duraj_random_2014} for asymptotic exit time formulas. When their results apply, the singular exponent is the very same irrational $\alpha$ or a simple rational function of it, and so the corresponding generating functions are not D-finite. A well-chosen weighting of the steps gives partial results when these theorems do not apply directly.  Central to the argument is the idea of the regime of the drift of the model, which we define next. 

\subsection{Drift}
\label{subsec:drift}
A key driver of the asymptotic behavior of both the counting sequence and the probability is the \emph{drift} of a model. The drift of a step set is the vector sum of its steps. If the steps are weighted (say, by a probability), then the drift is the appropriately weighted vector sum.  
Fix a step set~$S$ and let $\delta_S$ denote its drift.  The drift vector $\delta_S$ can be computed from the inventory of the model $P(x,y)$ by $\delta_S=\nabla(P)(1,1)=(P_x(1,1),P_y(1,1))$. In probabilistic terms, the drift is the expected direction of a step.

Following~\cite{Banuelos_Brownian_1997}, given a proper, open and connected subset $\Theta$ of the unit circle $\mathbb{S}^1\subset \mathbb{R}^2$, we can define the cone $C$ generated by $\Theta$ to be the set of all rays emanating from the origin and passing through $\Theta$:
\(
    C:= \{r\theta \mid r > 0, \theta \in \Theta\}.
\)
The \emph{polar cone} of $C$, denoted by $C^\#$, is defined by 
\begin{equation} \label{def:polarcone}
    C^\#:= \{x\in \mathbb{R}^2:\langle x, y\rangle \leq 0, \forall y\in C\},
\end{equation}
where $\langle x, y\rangle$ denotes the standard inner product. For a subset $A\subseteq \mathbb{R}^2$, we denote its interior by $A^{\circ}$. 

To each step set~$S$ we associate a cone  by means of a linear transformation $M$. The location of the drift vector relative to the cone associated to the model drives the dominant term of the coefficient asymptotics. If the model is not one of the five so-called singular models\footnote{None of the models considered here are singular.}, then the Laurent polynomial $P(x,y)$ has a unique positive critical point $(x_0, y_0)$ characterized by $P_x(x,y)=P_y(x,y)=0$. We define 
\begin{equation}\label{eq:c}
c:=\frac{P_{xy}(x_0,y_0)}{\sqrt{P_{xx}(x_0,y_0)P_{yy}(x_0,y_0)}},
\end{equation} and from it the angle $\theta=\arccos(-c)$.  A key result of \cite{bostan_non-d-finite_2014} that we use here is an exhaustive case analysis that proves the following about $c$.
\begin{theorem}[Bostan et al.~\cite{bostan_non-d-finite_2014}] \label{thm:irrational}
For each of the 51 non-singular infinite group 
models, the associated $c$ computed in Eq.~\eqref{eq:c} has the property that $-1-\frac{\pi}{\arccos(-c)}$ is irrational.
\end{theorem}
The matrix 
\[M:=\frac{1}{2 \sqrt{1-c^{2}}} 
\left[\begin{array}{cc}
\sqrt{1+c}+\sqrt{1-c} & \sqrt{1-c}-\sqrt{1+c} 
\\
 \sqrt{1-c}-\sqrt{1+c} & 2 \sqrt{1+c} 
\end{array}\right]
\]
is a linear transformation with the property that the model whose steps are the images of $S$ under $M$ has the identity as its covariance matrix. We define the cone \(C\) associated to a step set as the image of the first quadrant under \(M\), which is a wedge of angle \(\theta=\arccos(-c)\). The position of the drift vector relative to \(C\) and its polar cone \(C^\#\) determines the regime of the model.

\begin{example}
Consider model~9, with step set 
$\mathcal S=\{(-1,0),(0,-1),(0,1),(1,-1)\}$. Its inventory is
$P(x,y)=x^{-1}+y^{-1}+y+xy^{-1}$, and its drift is
$\delta_{\mathcal S}=(0,-1)$. The positive critical point is
$(x_0,y_0)=(r,r^2)$, where $r$ is the positive root of
$r^4-r-1=0$, giving $c=-1/(2r^{3/2})\approx -0.3707$ and
\[
M\approx
\begin{pmatrix}
1.0574 & 0.2032\\
0.2032 & 0.8541
\end{pmatrix}.
\]
Thus $C$ is generated by the column vectors
$v_1=(1.0574,0.2032)$ and $v_2=(0.2032,0.8541)$.  One can check that
$\langle\delta_{\mathcal S},v_1\rangle<0$ and
$\langle\delta_{\mathcal S},v_2\rangle<0$ so that
$\delta_{\mathcal S}\in(C^\#)^\circ$, and hence this model has a
polar interior drift regime. The six possible regimes are summarized
in Table~\ref{tab:regimes}.
\end{example}
\begin{table}
\centering
\begin{tabular}{lc@{\qquad}lc}
\toprule
Regime & Location of $\delta_S$
       & Regime & Location of $\delta_S$ \\
\midrule
interior drift
    & $C$
    & boundary drift
    & $\partial C\setminus\{0\}$ \\

non-polar exterior drift
    & $\mathbb{R}^{2}\setminus(C\cup C^\#)$
    & zero drift
    & $(0,0)$ \\

polar interior drift
    & $(C^\#)^\circ$
    & polar boundary drift
    & $\partial C^\#\setminus\{0\}$ \\
\bottomrule
\end{tabular}
\caption{Drift regimes}
\label{tab:regimes}
\end{table}

For models in two of these regimes, we appeal to known results on the exit time of random walks from cones to determine the singular exponent. 
\subsection{The probability model}
\label{subsec:theprobabilitymodel}
Consider a random walk $s(n)=X(1)+\dots+X(n)$, where $\{X(n)\}$ is a family of independent and identically distributed random variables taking values in some step set $S$. Let $\mathcal{K}=(0,\infty)^{2}$ be the open positive quadrant. For $x\in \mathcal{K}\cap \mathbb{Z}^{2}$, we denote by $\tau_x$ the exit time from $\mathcal{K}$ of the random walk starting at $x$: $\tau_x=\inf\{n\geq 1: x+s(n)\not\in \mathcal{K}\}$. Taking $e=(1,1)$, quarter-plane walks correspond precisely to random walks satisfying $\tau_e>n$. When the steps are chosen uniformly from $S$, we have that $\mathbf{P}(\tau_e>n) =\frac{c_S(n)}{|S|^n}$.
Thus, asymptotic estimates for $\mathbf P(\tau_e>n)$ give asymptotic
estimates for  $c_S(n)$ with the same singular exponent.

\subsection{Models with zero drift}
It is straightforward to iterate through all models and compute their drift. Doing so reveals that 
among the 56 infinite group cases, three have zero drift. They are listed in Table~\ref{tab:zerodrift}.
\begin{table}
\begin{tabular}{|M{0.8cm}|M{1.8cm}|M{2.5cm}|M{2.8cm}|M{1.5cm}|M{2cm}|} \hline
\(\#\) 
& Model 
& OEIS\\ \hline
30 &%
\StepModel{-1/-1,-1/0,0/1,1/-1,1/1} 
&\href{https://oeis.org/A151301}{\texttt{A151301}} \\
\hline
40 &
\StepModel{-1/-1,-1/1,0/1,1/0,1/-1} 
&\href{https://oeis.org/A151288}{\texttt{A151288}} \\
\hline
42 &
\StepModel{-1/0,-1/1,0/-1,1/1,1/-1} 
&\href{https://oeis.org/A151276}{\texttt{A151276}} \\
\hline
\end{tabular}
\caption{Zero drift models.}
\label{tab:zerodrift}
\end{table}
In each of these three cases, the conditions of Theorem 1 of Denisov and Wachtel~\cite[Equation 12]{denisov_random_2015} are met: the moment assumption holds (as the number of steps is finite), the normalization assumption holds since we can directly compute that the drift is \((0,0)\), and after the transformation $M$ the covariance matrix is the identity. The models are each aperiodic because they can reach every point. The cone is convex. 

In Example~2, Denisov and Wachtel give a formula for $\mathbf{P}(\tau_x>n)$. By their Theorem~1, \(\mathbf{P}(\tau_x>n) \sim \kappa V(x) n^{-p/2}\), for some constant $\kappa$ and a finite harmonic function $V(x)$. In this case $p$ is $\frac{\pi}{\arccos(-c)}$ with $c$ defined in Eq.~\eqref{eq:c}. The singular exponent and is irrational since it is half of $\frac{\pi}{\arccos(-c)}$, which is irrational (by Theorem~\ref{thm:irrational}). By Theorem~\ref{thm:alphaisirrational} we have shown the following result. 
\begin{theorem}
\label{thm:zerodrift}\label{thm:zerodriftnondfinite} 
The three small-step models with infinite group and zero drift each have a non-D-finite counting generating function. 
\end{theorem}

\subsection{Models with polar interior drift}
\label{subsec:polarintdrift}
For models with polar interior drift, we apply Theorem~1 of Duraj~\cite{duraj_random_2014} to give an asymptotic estimate for \(\mathbf{P}(\tau_x>n)\). The applicability is clarified in~\cite[Example 7]{duraj_random_2014}.  

For random walks with nonzero drift, we use the following specialization of Duraj's theorem~\cite{duraj_random_2014} to the cone $\mathcal{K}$, for which the additional geometric assumptions are satisfied. 
\begin{theorem}[Duraj 2014, Theorem 1~{\cite{duraj_random_2014}}] Suppose that a random walk has nonzero drift and let $R(h):=\mathbb E\!\left[e^{\langle h,X\rangle}\right]$ be the moment generating function of the random step \(X\), where $h\in\mathbb R^2$. Assume that the following conditions hold:  
\begin{description}
\item[\rm(A1) Cramér condition]
There exists a nonzero vector $h \in \mathcal{K}$ such that $\nabla(R(h))=0$,
\item[\rm(A2)  non-collinearity]
The step set is not contained in a line passing through the origin. Equivalently, $\mathbf P(\langle u,X\rangle=0)<1$ for every nonzero \(u\in\mathbb R^2\),
\item[\rm(A3)  strong aperiodicity]
    For every \(x\in\mathbb Z^2\), the smallest additive subgroup of \(\mathbb Z^2\) containing  the set \(\{x+z:\mathbf P(X=z)>0\}\) is \(\mathbb Z^2\).
\end{description}
Then, for each fixed endpoint set \(A\subset \mathcal{K}\cap\mathbb Z^2\), there exists a constant \(C_A>0\) such that
\[
\mathbf P\bigl(e+s(n)\in A,\ \tau_e>n\bigr)
\sim
C_A\,R(h)^n n^{-(p+1)},
\qquad n\to\infty.
\]
Moreover, $p=\pi/\theta$, where $\theta$ is the opening angle of the covariance-normalized cone; for the cone \(\mathcal{K}\), \(\theta=\arccos(-c)\).

\label{thm:duraj}
\end{theorem}
\begin{theorem}
\label{thm:interiordriftnondfinite}
The eleven small-step models in Table~\ref{tab:polardrift}, with infinite group and polar interior drift have a non D-finite counting generating function. 
\end{theorem}
\begin{proof}
Duraj's result applies to models with drift regime in $(C^\#)^\circ$, and the main work of this proof is to identify these cases. Again, it is straightforward linear algebra to iterate through all cases, compute the drift regime for each model, and output the figure of the cone and the drift.  The Cram\'er condition is easily verified computationally. In this case we need to verify that $x_0>1, y_0>1$. 
Table~\ref{tab:polar-critical} provides approximate values, and they are all greater than 1.

\begin{table}[ht]
\centering
\begin{tabular}{ccc@{\qquad}ccc}
\toprule
Model $m$ & $x_0$ & $y_0$ & Model $m$ & $x_0$ & $y_0$ \\
\midrule
9  & 1.220744 & 1.490216 & 36 & 1.985438 & 1.286606 \\
10 & 1.577018 & 1.486986 & 47 & 1.244924 & 1.065750 \\
16 & 1.220744 & 1.220744 & 48 & 2.278288 & 1.607273 \\
24 & 1.324718 & 1.324718 & 53 & 1.220744 & 1.220744 \\
33 & 1.676023 & 1.809053 & 55 & 1.324718 & 1.324718 \\
35 & 1.754878 & 1.324718 &    &          &          \\
\bottomrule
\end{tabular}
\caption{Approximate positive critical points $(x_0,y_0)$ of the inventory for
the polar-interior models of Theorem~11. In each case $x_0>1$ and $y_0>1$.}
\label{tab:polar-critical}
\end{table}

The non-collinearity is trivially satisfied for all non-trivial models. 
We certify the aperiodicity in Table~{tab:polar-aperiodic}

\begin{table}[ht]
\centering
\begin{tabular}{ccc}
\toprule
Model(s) $m$ & $u$ & $v$ \\
\midrule
9                  & $\leftarrow-\uparrow=(-1,-1)$   & $\downarrow-\searrow=(-1,0)$ \\
10                 & $\leftarrow-\uparrow=(-1,-1)$   & $\leftarrow-\swarrow=(0,1)$ \\
16, 33, 48, 53, 55 & $\swarrow-\leftarrow=(0,-1)$    & $\swarrow-\downarrow=(-1,0)$ \\
24                 & $\leftarrow-\swarrow=(0,1)$     & $\downarrow-\swarrow=(1,0)$ \\
35                 & $\leftarrow-\nwarrow=(0,-1)$    & $\nwarrow-\uparrow=(-1,0)$ \\
36, 47             & $\swarrow-\leftarrow=(0,-1)$    & $\nwarrow-\uparrow=(-1,0)$ \\
\bottomrule
\end{tabular}
\caption{Step differences $u$ and $v$ certifying strong aperiodicity for the
eleven polar-interior models of Theorem~11. In each row, $|\det(u,v)|=1$.}
\label{tab:polar-aperiodic}
\end{table}

By Duraj's theorem we find that the singular exponent for every model in these cases is $-1-\pi/\arccos(-c)$ (for the respective value of $c$), which is irrational by Theorem~\ref{thm:irrational}. 
The result follows by Theorem~\ref{thm:alphaisirrational}.
\end{proof}
 
\begin{longtable}{|M{0.8cm}|M{1.8cm}|M{1.5cm}|M{2.2cm}|M{2.5cm}|M{2.2cm}|} 
\hline
\(\#\) 
& Model 
& Drift Regime
& OEIS \\
\hline
\endfirsthead

\hline
\multicolumn{6}{|c|}{\tablename~\thetable{} -- continued} \\
\hline
\(\#\) 
& Model 
& Drift Regime

& OEIS \\
\hline
\endhead

\hline
\multicolumn{6}{|r|}{Continued on next page} \\
\hline
\endfoot

\endlastfoot

9 
&\StepModel{-1/0,0/1,0/-1,1/-1} 
&\DriftAndRegion{-1/0,0/1,0/-1,1/-1}{(1.0,0.15)}{(.16,0.98)}{(0.15,-1)}{(-.98,.16)}
&\href{https://oeis.org/A151262}{\texttt{A151262}} \\
\hline

10 &
\StepModel{-1/0,0/1,-1/-1,1/-1} & \DriftAndRegion{-1/0,0/1,-1/-1,1/-1}{(1.,.12)}{(.14,1.0)}{(.12,-1.)}{(-1.0,.14)}
&\href{https://oeis.org/A151257}{\texttt{A151257}} \\
\hline


16 &
\StepModel{-1/-1,-1/0,0/-1,1/1} & \DriftAndRegion{-1/-1,-1/0,0/-1,1/1}{(1.3,-.51)}{(-.27,.95)}{(-.51,-1.3)}{(-.95,-.27)}
&\href{https://oeis.org/A151270}{\texttt{A151270}} \\
\hline
24 &
\StepModel{-1/0,0/1,0/-1,1/0,-1/-1} & \DriftAndRegion{-1/0,0/1,0/-1,1/0,-1/-1}{(1.,-.11)}{(-.10,1.0)}{(-.12,-1.)}{(-1.0,-.10)}
&\href{https://oeis.org/A151286}{\texttt{A151286}}\\
\hline
33 &
\StepModel{-1/-1,-1/0,0/1,0/-1,1/-1} & \DriftAndRegion{-1/-1,-1/0,0/1,0/-1,1/-1}{(1.,.12)}{(.14,1.0)}{(.12,-1.)}{(-1.0,.14)}
&\href{https://oeis.org/A151263}{\texttt{A151263}} \\
\hline


35 &
\StepModel{-1/0,-1/1,0/1,0/-1,1/-1} & \DriftAndRegion{-1/0,-1/1,0/1,0/-1,1/-1}{(1.2,.45)}{(.50,.86)}{(.45,-1.2)}{(-.86,.50)}
&\href{https://oeis.org/A151268}{\texttt{A151268}} \\
\hline

36 &
\StepModel{-1/-1,-1/0,-1/1,0/1,1/-1} & \DriftAndRegion{-1/-1,-1/0,-1/1,0/1,1/-1}{(1.1,.31)}{(.34,.92)}{(.31,-1.1)}{(-.92,.34)}
&\href{https://oeis.org/A151260}{\texttt{A151260}} \\
\hline

47 &
\StepModel{-1/-1,-1/0,-1/1,0/1,1/0,1/-1} & \DriftAndRegion{-1/-1,-1/0,-1/1,0/1,1/0,1/-1}{(1.,.12)}{(.14,1.0)}{(.12,-1.)}{(-1.0,.14)}
&\href{https://oeis.org/A151298}{\texttt{A151298}} \\
\hline
48 &
\StepModel{-1/-1,-1/0,-1/1,0/1,0/-1,1/-1} & \DriftAndRegion{-1/-1,-1/0,-1/1,0/1,0/-1,1/-1}{(1.12,-.30)}{(-.21,.97)}{(-0.31,-1.12)}{(-.97,-.21)} 
&\href{https://oeis.org/A151269}{\texttt{A151269}} \\
\hline

53 &
\StepModel{-1/-1,-1/0,-1/1,0/-1,1/1,1/-1} & \DriftAndRegion{-1/-1,-1/0,-1/1,0/-1,1/1,1/-1}{(1.,0)}{(0,1.)}{(0,-1.)}{(-1.,0)}
&\href{https://oeis.org/A151277}{\texttt{A151277}} \\
\hline

55 &\raisebox{5pt}[0pt][0pt]{\StepModel{-1/-1,-1/0,-1/1,0/1,0/-1,1/-1,1/0}} & \DriftAndRegion{-1/-1,-1/0,-1/1,0/1,0/-1,1/-1,1/0}{(1.0,.14)}{(.15,1.0)}{(.14,-1.0)}{(-1.0,.15)}

&\href{https://oeis.org/A151305}{\texttt{A151305}} \\
\hline

\caption{Models with polar-interior drift with their cone}
\label{tab:polardrift}
\end{longtable}

\section{Some non-polar exterior and interior-drift models with a non-D-finite boundary series}
\label{sec:bdrynondfinite}

We now turn to a further family of 18 infinite group models for which the preceding arguments do not establish the non-D-finiteness of the total counting series $\Qtot$, since they have neither zero nor polar interior drift. However, by applying a special weighting we can satisfy the conditions of Duraj's theorem applied to walks killed on a boundary. Each of these models has drift that is zero in one component and is sufficiently small  in the other so that using a well chosen weighting we can connect it to a model that does have polar interior drift in a very direct way. It is closely related to the strategy in~\cite{bostan_non-d-finite_2014}.

Unlike the asymptotic estimates for survival probabilities alone, Duraj's theorem allows one to impose restrictions on the endpoint of the walk. This provides a way to study the asymptotic behavior of different counting sequences according to their endpoint conditions, rather than only the total number of walks confined to the quarter plane. 

For each model in this case, we show that the asymptotics of $\Qx$ involve, yet again, the same irrational number as in Theorem~\ref{thm:irrational}, which proves that~$\Qx$ is not D-finite.  To do this, we next define a specific weighting on each step to create a model with polar interior drift. The weighted model is chosen to preserve $\Qx$, allowing us to obtain the required coefficient asymptotics.

\subsection{Weighted models as random walks}
\label{subsec: weightedrandomwalks}
We use central weightings in a similar manner to~\cite{courtiel_weighted_2017}. A central weighting assigns to each step \(s=(i,j)\in S\) a positive weight of the form $wt(s)=a^i b^j$, where $a,b>0.$ The weight $wt(\omega)$ of an \(n\)-step walk \(\omega=(\omega_1,\ldots,\omega_n)\) is the product of the weights of its steps. If \(\omega\) ends at \((u,v)\), then its weight is $wt(\omega):=\prod_{k=1}^n wt(\omega_k)=a^u b^v$, which depends only on its endpoint.

Let \(Q^{w}(x,y;t)\) denote the endpoint generating function of the corresponding weighted quarter-plane model. Since every walk ending at \((u,v
)\) gets the factor \(a^u b^v\), we have the endpoint generating function identity
\begin{equation}
\label{eq:weigthedendpointgf}
Q^w(x,y;t)=Q(ax,by;t).    
\end{equation}
Specializing Eq.~\eqref{eq:weigthedendpointgf} at $x=1$ and $y=0$, we obtain $$Q^{w}(1,0;t)=\left.Q(ax,by;t)\right|_{x=1,y=0}=Q(a,0;t).$$ If we set $wt((i,j))=b^j$ for some positive $b$, then $Q^{w}(1,0;t)=Q(1,0;t)$ as series, and so their coefficients are equal.   We can convert this into a probability distribution by normalizing the weights by $W$, defined by
\[
W\:= \sum_{s\in S}wt(s)=P(1,b),
\qquad
p_s \:= \frac{wt(s)}{W},
\qquad s\in\mathcal S.
\]

We equip the random walk introduced in Section~\ref{subsec:theprobabilitymodel} with step probabilities $p_s$. Since the random increments $X(k)$ are independent, the probability of an $n$-step quarter-plane walk \(\omega=(\omega_1,\ldots,\omega_n)\) is
\begin{equation}
\mathbf P\bigl(X(1)=\omega_1,\ldots,X(n)=\omega_n\bigr)
=
\prod_{k=1}^n p_{\omega_k}
=
\frac{wt(\omega)}{W^n}.    
\label{eq:weightmodelprobability}
\end{equation}

In particular, every $n$-step quarter-plane walk ending on the horizontal axis $B=\{(i,0): i\geq 0 \}$ has weight $1$ and hence probability $W^{-n}$. Since $[t^{n}]Q(1,0;t)$ counts precisely these walks, summing their probabilities yields
\begin{equation}
\label{eq: weightedmodelcounting}
[t^{n}]Q(1,0;t)=W^n\mathbf{P}\left(e+s(n)\in e+B, \tau_{e}>n \right).
\end{equation}

\subsection{Non-D-finiteness of the boundary series}
 We can apply this idea to establish the non-$D$-finiteness of the boundary generating functions $Q_{S_{m}}(1,0;t)$, for the models with non-polar exterior drift listed in Table~\ref{tab:exterior} and those with interior drift listed in Table~\ref{tab:intnondfinite}, with the exception of models  \(7,17,27,\) and \(32\) because they fail the strong aperiodicity condition\footnote{Aperiodicity can often be managed by considering instead the step set $S+S$ and looking at walks of length $2n$ rather than $n$, as in~\cite{bostan_non-d-finite_2014}. } of Theorem~\ref{thm:duraj}. Duraj~\cite{duraj_random_2014} shows that this assumption can be relaxed for exit-time asymptotics, but without the endpoint restrictions needed here.

\begin{theorem}
\label{thm:boundary-series-nondfinite}
Let $S_m$ denote the step set of model $m$. For each of the eighteen strongly aperiodic models $m$ listed in Tables~\ref{tab:exterior} and \ref{tab:intnondfinite}, the
boundary generating function $Q_{S_m}(1,0;t)$ is not D-finite.
\end{theorem}

\begin{proof}
Fix a model $m$ among those considered in the theorem, with step set $S_m$, inventory $P_{S_m}$ and positive critical point $(x_m,y_m)$. As established in Section~\ref{subsec: weightedrandomwalks}, a central weighting $wt_{m}{(i,j)}=b_{m}^{j}$, where $b_m>0$ will be chosen below, preserves the boundary generating function $Q_{S_m}(1,0;t)$. We consider the random walk associated with this weighting, whose step distribution is 
$$\mathbf{P}(X=(i,j))=\frac{b_m^{j}}{P_{S_m}(1,b_m)}, \ (i,j)\in S_m.$$

We now verify that this random walk satisfies conditions {\rm (A1)--(A3)} of Theorem~\ref{thm:duraj}.\\
\emph{Verification of~{\rm (A1)}}. For $h=(h_1,h_2)\in \mathbb{R}^{2}$, the moment generating function of the associated random walk is \[R_m(h_1,h_2) = \frac{P_{S_{m}}(e^{h_1},b_me^{h_2})}{P_{S_{m}}(1,b_m)}. \] Since $(x_m,y_m)$ is a critical point of $P_{S_m}$, the vector \[ h_m:= \left( \log x_m, \log \frac{y_m}{b_m} \right) \] satisfies $\nabla R_m(h_m)=0$. For each model $m$ under consideration, \(x_m>1\). Hence, choosing $0<b_m<y_m$ ensures that $h_m\in \mathcal{K}=(0,\infty)^{2}$, as required by the Cramér condition. The approximate values of $y_m$ are recorded in Table~\ref{tab:admissible-weights}; in particular, \(b_m=1/3\) is admissible for every model in Theorem~\ref{thm:boundary-series-nondfinite}.

\begin{table}[ht]
\centering
\begin{tabular}{ccc@{\qquad}ccc}
\toprule
Model $m$ & $x_m$ & $y_m$ & Model $m$ & $x_m$ & $y_m$ \\
\midrule
3  & 1.490216 & 0.819173 & 26 & 1.754878 & 0.754878 \\
4  & 1.486986 & 0.634108 & 28 & 1.676023 & 0.552775 \\
6  & 1.220744 & 0.671044 & 37 & 1.809053 & 0.596651 \\
12 & 1.220744 & 0.819173 & 41 & 1.985438 & 0.777238 \\
18 & 1.577018 & 0.672501 & 43 & 1.607273 & 0.438926 \\
20 & 1.286606 & 0.503667 & 44 & 1.065750 & 0.803262 \\
23 & 1.324718 & 0.754878 & 49 & 2.278288 & 0.622172 \\
25 & 1.324718 & 0.569840 & 50 & 1.244924 & 0.938306 \\
   &          &          & 52 & 1.220744 & 0.819173 \\
   &          &          & 54 & 1.324718 & 0.754878 \\
\bottomrule
\end{tabular}
\caption{Approximate positive critical points $(x_m,y_m)$ for the models of
Theorem~12. In each case $x_m>1$, and the admissible range for the weight is
$0<b_m<y_m$; in particular $b_m=1/3$ is admissible for every model.}
\label{tab:admissible-weights}
\end{table}


\emph{Verification of {\rm (A2)}.}   
Since $b_m>0$, every step in \(\mathcal S_m\) has positive probability. Thus, for any nonzero vector $u\in \mathbb{R}^{2}$,
\(\mathbf P(\langle u,X\rangle=0)=1\) if and only if \(
\langle u,s\rangle=0\) for every \(s\in S_m.\) Consequently, condition {\rm (A2)} is equivalent to $\operatorname{span}_{\mathbb{R}}(S_m)=\mathbb{R}^{2}$. Each model under consideration contains two linearly independent steps, so {\rm (A2)} holds.

\emph{Verification of {\rm (A3)}.} For any
$x\in\mathbb Z^2$ and $s,s' \in S_m$, we have that $(x+s)-(x+s')=s-s'$. Thus, the group generated by $\{(x+s):s\in S_m\}$ contains all step differences. Consequently, it suffices to show that $\left\langle s-s':s,s'\in S_m\right\rangle_{\mathbb Z}=\mathbb Z^2.$
For each model in Theorem~\ref{thm:boundary-series-nondfinite}, Table~\ref{tab:aperiodicity} gives two step differences \(u,v\) satisfying $|\det(u,v)|=1.$ Thus \(u\) and \(v\) form a \(\mathbb Z\)-basis of \(\mathbb Z^2\), and condition {\rm (A3)} follows.
\begin{table}[ht]
\centering
\renewcommand{\arraystretch}{1.2}
\begin{tabular}{ccc}
\toprule
Model(s) \(m\) & \(u\) & \(v\)\\
\midrule
\(3\)
& \(\rightarrow-\uparrow=(1,-1)\)
& \(\downarrow-\swarrow=(1,0)\)\\

\(4\)
& \(\rightarrow-\uparrow=(1,-1)\)
& \(\rightarrow-\searrow=(0,1)\)\\

\(6\)
& \(\uparrow-\nearrow=(-1,0)\)
& \(\uparrow-\leftarrow=(1,1)\)\\

\(12\)
& \(\uparrow-\nwarrow=(1,0)\)
& \(\nwarrow-\leftarrow=(0,1)\)\\

\(18\)
& \(\nwarrow-\leftarrow=(0,1)\)
& \(\downarrow-\leftarrow=(1,-1)\)\\

\(20,25,43,44,54\)
& \(\rightarrow-\nearrow=(0,-1)\)
& \(\uparrow-\nearrow=(-1,0)\)\\

\(23,37\)
& \(\rightarrow-\searrow=(0,1)\)
& \(\downarrow-\searrow=(-1,0)\)\\

26
& \(\uparrow-\nearrow=(-1,0)\)
& \(\swarrow-\leftarrow=(0,-1)\)
\\

\(28,49,52\)
& \(\uparrow-\nearrow=(-1,0)\)
& \(\nwarrow-\leftarrow=(0,1)\)\\

\(41\)
& \(\nwarrow-\leftarrow=(0,1)\)
& \(\downarrow-\swarrow=(1,0)\)\\

50
& \(\rightarrow-\nearrow=(0,-1)\)
& \(\downarrow-\swarrow=(1,0)\)
\\
\bottomrule
\end{tabular}

\caption{Step differences \(u\) and \(v\) certifying strong
aperiodicity. In each row, \(\lvert\det(u,v)\rvert=1\)}
\label{tab:aperiodicity}
\end{table}

Having verified {\rm (A1)}--{\rm (A3)}, we may apply Theorem~\ref{thm:duraj}
with $e=(1,1)$ and $A=e+B=\{(i,1), i\geq 1\}$. The covariance-normalized cone has opening angle
\[
\theta_m=\arccos(-c_m),
\qquad
c_m\coloneqq
\frac{P_{{S_m},xy}(x_m,y_m)}
{\sqrt{P_{{S_m},xx}(x_m,y_m)P_{{S_m},yy}(x_m,y_m)}}.
\]

Moreover, $P_{S_m}(1,b_m)R_m(h_m)=P_{S_m}(x_m,y_m)$. Consequently, Theorem~\ref{thm:duraj} and Eq.~\eqref{eq: weightedmodelcounting} give
$$[t^n]Q_{S_m}(1,0;t) \sim C_m\rho_m^n n^{-(p_m+1)},$$
for some $C_m>0,$ where $\rho_m=P_{S_m}(x_m,y_m)$ and $p_m=\pi/\theta_m$. The singular exponent $p_m$ also appears in the excursion asymptotics of Bostan, Raschel and Salvy~\cite{bostan_non-d-finite_2014}. By Theorem~\ref{thm:irrational}, $p_m$ is irrational and hence so is $-(p_m+1)$. Therefore, by Theorem~\ref{thm:alphaisirrational} $Q_{S_m}(1,0;t)$ is not D-finite.
\end{proof}

\section{Conjectured singular exponents}
When the number of walks has the asymptotic form $\kappa \rho^{-n} n^\alpha$, we can deduce $\alpha$. This is true for a number of regimes but has not been proven in general. The results for zero drift~\cite{denisov_random_2015}, polar interior~\cite{duraj_random_2014} and some of interior ~\cite{nguyen_random_2026} regimes are proven, as are the $\alpha_{00}$ for all cases.  We have strong computational evidence for this table, under the assumption of this asymptotic shape. Note that the case distinctions follow those of Fayolle and Raschel~\cite{fayolle_some_2012} for the dominant singularity. 
\begin{table}[ht]
\centering
\renewcommand{\arraystretch}{1.15}
\begin{tabular}{llccc cccc}
\toprule
\multirow{2}{*}{Drift regime}
& \multirow{2}{*}{\(\#\)}
& \multicolumn{3}{c}{Sign conditions}
& \multicolumn{4}{c}{Coefficient exponents}
\\
\cmidrule(lr){3-5}
\cmidrule(lr){6-9}
&
& \(\delta_{S,x}\)
& \(\delta_{S,y}\)
& \(\kappa_S\)
& \(-\alpha_{11}\)
& \(-\alpha_{10}\)
& \(-\alpha_{01}\)
& \(-\alpha_{00}\)
\\
\midrule

\multirow{2}{*}{Interior}
& \multirow{2}{*}{13}
& \(+\) & \(+\) & $\ast$
& \(0\) & \(3/2\) & \(3/2\) & \(p+1\)
\\
&
& \(0\) & \(+\) & \(+\)
& \(1/2\) & \(p+1\) & \(3/2\) & \(p+1\)
\\

\midrule

\multirow{4}{*}{Non-polar exterior}
& \multirow{4}{*}{22}
& \(-\) & \(+\) & $\ast$
& \(3/2\) & \(p+1\) & \(3/2\) & \(p+1\)
\\
&
& \(-\) & \(0\) & \(+\)
& \(3/2\) & \(p+1\) & \(3/2\) & \(p+1\)
\\
&
& \(+\) & \(0\) & \(-\)
& \(1/2\) & \(3/2\) & \(3/2\) & \(p+1\)
\\
&
& \(0\) & \(+\) & \(-\)
& \(1/2\) & \(3/2\) & \(3/2\) & \(p+1\)
\\

\midrule

Zero drift
& 3
& \(0\) & \(0\) & $\ast$
& \(p/2\) & \(1+p/2\) & \(1+p/2\) & \(p+1\)
\\

\midrule

\multirow{3}{*}{Polar interior}
& \multirow{3}{*}{13}
& \(-\) & \(-\) & $\ast$
& \(p+1\) & \(p+1\) & \(p+1\) & \(p+1\)
\\
&
& \(0\) & \(-\) & \(-\)
& \(p+1\) & \(p+1\) & \(p+1\) & \(p+1\)
\\
&
& \(-\) & \(0\) & \(-\)
& \(p+1\) & \(p+1\) & \(p+1\) & \(p+1\)
\\

\bottomrule
\end{tabular}

\caption{Conjectured critical exponents for the 51 nonsingular infinite-group
models, where \(p=\frac{\pi}{\arccos(-c_S)}.\) Here \(\kappa_S=\operatorname{Cov}(X_1,X_2)\), where
\(X=(X_1,X_2)\) is a uniformly chosen step in \(S\). An $\ast$ indicates that the sign of \(\kappa_S\) is irrelevant. Entries for zero, polar interior drift and $\alpha_{00}$ are proven}
\label{tab:singular-exponents}
\end{table}

When exactly one drift coordinate vanishes and the other is positive, the sign of \(\kappa_S\) distinguishes the interior cases in Table~9 from the non-polar exterior cases in Table~10. This classification also identifies the limitation of the present methods. For zero-drift and polar-interior models, the irrational exponent occurs in the asymptotics of the coefficients of $Q(1,1;t)$, and therefore prove its non-$D$-finiteness. For the models with non-polar exterior drift and interior drift, the irrational exponent occurs instead in the asymptotics of the coefficients of a boundary generating function. Although this proves this series is non-$D$-finite, passing from this result to the total counting series remains open.

\section{Conclusion and open questions}
\label{sec:conclusion}
Our results establish the non-D-finiteness of $Q(1,1;t)$ for infinite-group models with zero or polar-interior drift. For certain models with non-polar exterior or interior drift, we instead establish the non-D-finiteness of the generating function $Q(1,0;t)$. A natural direction is to develop methods that transfer the non-D-finiteness obstructions found in $Q(1,0;t)$ to $Q(1,1;t)$. More generally, resolving the remaining cases of Conjecture~\ref{conj:q11} may require new techniques, perhaps adapting the differential transcendence arguments of~\cite{hardouin_d-finiteness_2025}. Six of the nine decoupled, differentially algebraic models (3, 4, 6, 25, 26, 37) are among the models in Theorem 12. 




\paragraph{Author contributions}

\begin{itemize}
\item MM: supervision, project administration, conceptualization, Writing: original draft, editing.
\item LP: conceptualization, proof arguments, validation, writing: review and editing.
\item JP: conceptualization, proofs, validation, writing: original draft, review and editing. 
\end{itemize}
We gratefully acknowledge our funders: Canada NSERC Discovery Grant "Combinatorics and Transcendence" and an SFU Dean of Science Research Grant and a CNRS MITI grant.
\printbibliography
\newpage
\appendix
\section{Singular exponents}
In the tables that follow, $\alpha_{i,j}$ denotes the singular exponent associated with $Q(i,j;t)$, assuming an asymptotic form $[t^n]Q(i,j;t) \sim K \rho^n n^{\alpha_{i,j}}$ for constants $K$ and $\rho$. Models that fail the strong aperiodicity condition (A3) are marked with an asterisk. The tables are organized by the drift regime and value of $\alpha_{11}$.

\begin{longtable}{|M{0.8cm}|M{1.5cm}|M{1.5cm}|M{1.5cm}|M{1.5cm}|M{1.5cm}|M{1.5cm}|M{2cm}|} 
\hline
\(\#\) 
& Model 
& Drift
& $-\alpha_{11}$
& $-{\alpha_{10}}$
& $-{\alpha_{01}}$
& $-{\alpha_{00}}$
& OEIS\\
\hline
\endfirsthead

\hline
\multicolumn{8}{|c|}{\tablename~\thetable{} -- continued} \\
\hline
\(\#\) 
& Model 
& Drift 
& $-{\alpha_{11}}$
& $-{\alpha_{10}}$
& $-{\alpha_{01}}$
& $-{\alpha_{00}}$
& OEIS\\
\hline
\endhead

\hline
\multicolumn{8}{|r|}{Continued on next page} \\
\hline
\endfoot

\endlastfoot

3 &
\StepModel{-1/-1,1/0,-1/0,0/1} & \DriftAndRegion{-1/-1,1/0,-1/0,0/1}{(1.1,-.22)}{(-.17,1.0)}{(-.22,-1.1)}{(-1.0,-.17)}
&$3/2$ 
&2.6106
& $3/2$
&2.6106
&\href{https://oeis.org/A151279}{\texttt{A151279}} \\
\hline

4 &
\StepModel{-1/-1,1/0,0/1,-1/1} & \DriftAndRegion{-1/-1,1/0,0/1,-1/1}{(.99,-.11)}{(-.10,1.0)}{(-.11,-.99)}{(-1.0,-.10)}
& {$3/2$}
&2.7204
& {$3/2$}
&2.7204
& \href{https://oeis.org/A151280}{\texttt{A151280}} \\
\hline

12 &
\StepModel{-1/1,-1/0,0/1,1/-1} & \DriftAndRegion{-1/1,-1/0,0/1,1/-1}{(1.3,.53)}{(.58,.85)}{(.53,-1.3)}{(-.85,.58)}
&{$3/2$}
&{5.1362}
&$3/2$
&5.1362
&\href{https://oeis.org/A151259}{\texttt{A151259}} \\
\hline

17* &
\StepModel{-1/1,-1/-1,0/-1,1/1} & \DriftAndRegion{-1/1,-1/-1,0/-1,1/1}{(1.0,-.26)}{(-.20,1.0)}{(-.26,-1.0)}{(-1.0,-.20)}
& $3/2$
&{2.5193}
& $3/2$
& 2.5211
&\href{https://oeis.org/A151271}{\texttt{A151271}} \\
\hline

18 &
\StepModel{-1/0,-1/1,0/-1,1/1} & \DriftAndRegion{-1/0,-1/1,0/-1,1/1}{(.97,-.13)}{(-.12,1.0)}{(-.13,-.97)}{(-1.0,-.12)}
&$3/2$
&2.7204
&$3/2$
&2.7204
&\href{https://oeis.org/A151272}{\texttt{A151272}} \\
\hline

23 &
\StepModel{0/-1,1/0,-1/0,0/1,-1/1} & \DriftAndRegion{0/-1,1/0,-1/0,0/1,-1/1}{(1.01,-.109)}{(-.0973,1.00)}{(-.109,-1.01)}{(-1.00,-.0973)}
& $3/2$
& 3.3202
& $3/2$
&3.3202
&\href{https://oeis.org/A151295}{\texttt{A151295}} \\
\hline

26 &
\StepModel{0/1,1/1,-1/0,-1/-1,0/-1}&
\DriftAndRegion{0/1,1/1,-1/0,-1/-1,0/-1}
{(1.21,-.43)}{(-.43,1.21)}
{(-.43,-1.21)}{(-1.21,-.43)}
&
$3/2$
&
2.3976
&
$3/2$
&
2.3976
&
\href{https://oeis.org/A151299}{\texttt{A151299}}
\\
\hline

27* &
\StepModel{-1/-1,-1/1,0/1,0/-1,1/1} & \DriftAndRegion{-1/-1,-1/1,0/1,0/-1,1/1}{(1.0,-.27)}{(-.21,1.0)}{(-.27,-1.0)}{(-1.0,-.21)}
&$3/2$
&2.5027
& $3/2$
&2.5035
&\href{https://oeis.org/A151300}{\texttt{A151300}} \\
\hline

28 &
\StepModel{-1/1,-1/0,0/1,0/-1,1/1} & \DriftAndRegion{-1/1,-1/0,0/1,0/-1,1/1}{(.99,-.11)}{(-.10,1.0)}{(-.11,-.99)}{(-1.0,-.10)} 
& $3/2$
& 2.7421
& $3/2$
& 2.7421
& \href{https://oeis.org/A151303}{\texttt{A151303}} \\
\hline

32* &
\StepModel{-1/-1,1/0,-1/0,-1/1,1/1} & \DriftAndRegion{-1/-1,1/0,-1/0,-1/1,1/1}{(1.1,-.31)}{(-.22,.99)}{(-.31,-1.1)}{(-.99,-.22)} 
&$3/2$
&3.1595
&$3/2$
&2.5035
&\href{https://oeis.org/A151306}{\texttt{A151306}} \\
\hline

37 &
\StepModel{-1/-1,1/0,-1/0,0/1,-1/1} & \DriftAndRegion{-1/-1,1/0,-1/0,0/1,-1/1}{(.99,-.11)}{(-.10,1.0)}{(-.11,-.99)}{(-1.0,-.10)} 
&$3/2$
&2.7421
& $3/2$
&2.7421
& \href{https://oeis.org/A151285}{\texttt{A151285}} \\
\hline

41 &
\StepModel{-1/-1,-1/0,-1/1,0/-1,1/1} & \DriftAndRegion{-1/-1,-1/0,-1/1,0/-1,1/1}{(1.1,-.31)}{(-.22,.99)}{(-.31,-1.1)}{(-.99,-.22)} 
&$3/2$
&2.4829
&$3/2$
&2.4829
&\href{https://oeis.org/A151274}{\texttt{A151274}} \\
\hline

43 &
\StepModel{-1/-1,1/0,-1/0,1/1,0/1,-1/1}
&
\DriftAndRegion{-1/-1,1/0,-1/0,1/1,0/1,-1/1}
{(-.27,1.0)}{(1.0,-.21)}{(-1.0,-.27)}{(-.21,-1.0)}
&
$3/2$
&
2.4911
&
$3/2$
&
2.4911
&
\href{https://oeis.org/A151322}{\texttt{A151322}} \\
\hline

49 &
\StepModel{-1/-1,-1/0,-1/1,0/1,0/-1,1/1} & \DriftAndRegion{-1/-1,-1/0,-1/1,0/1,0/-1,1/1}{(.965,-.265)}{(-.211,.980)}{(-.265,-.965)}{(-.980,-.211)}
& $3/2$
& {2.4911}
&$3/2$
&2.4911
&\href{https://oeis.org/A151309}{\texttt{A151309}} \\
\hline

50 &
\StepModel{-1/-1,0/-1,1/0,-1/0,1/1,-1/1} & \DriftAndRegion{-1/-1,0/-1,1/0,-1/0,1/1,-1/1}{(.990,-.150)}{(-.131,.991)}{(-.150,-.990)}{(-.991,-.131)}
& $3/2$
& 2.6798
& $3/2$
& 2.6798
&\href{https://oeis.org/A151311}{\texttt{A151311}} \\
\hline

52 &
\StepModel{-1/-1,-1/0,-1/1,0/1,1/1,1/-1} & \DriftAndRegion{-1/-1,-1/0,-1/1,0/1,1/1,1/-1}{(1.,.0)}{(0,1.)}{(0,-1.)}{(-1.,0)} 
& $3/2$
& 3.0421
&$3/2$
& 3.0421
&\href{https://oeis.org/A151310}{\texttt{A151310}} \\
\hline

54 &
\StepModel{-1/-1,-1/0,-1/1,0/1,0/-1,1/1,1/0} & \DriftAndRegion{-1/-1,-1/0,-1/1,0/1,0/-1,1/1,1/0}{(.991,-.156)}{(-.134,.991)}{(-.156,-.991)}{(-.991,-.134)}
& $3/2$
& 2.6680
&$3/2$
&2.6680
&\href{https://oeis.org/A151328}{\texttt{A151328}} \\
\hline

\caption{The 17 non-polar exterior-drift models with ${-\alpha_{11}}\sim 1.5$ and ${-\alpha_{10}}$ irrational}

\label{tab:exterior}
\end{longtable}
 
\begin{longtable}{|M{0.8cm}|M{1.5cm}|M{1.5cm}|M{1.5cm}|M{1.5cm}|M{1.5cm}|M{1.5cm}|M{2cm}|} 
\hline
\(\#\) 
& Model 
& Drift
& $-{\alpha_{11}}$
& $-{\alpha_{10}}$
& $-{\alpha_{01}}$
& $-{\alpha_{00}}$
& OEIS\\
\hline

$6$ &
\StepModel{-1/0,0/1,0/-1,1/1} &      \DriftAndRegion{-1/0,0/1,0/-1,1/1}
{(.982026,-.188747)}
{(-.159161,.987253)}
{(-.188747,-.982026)}
{(-.987253,-.159161)}
&$1/2$
&2.6106
&$3/2$
&2.6106
&\href{https://oeis.org/A151290}{\texttt{A151290}} \\
\hline

7* &
\StepModel{-1/-1,1/0,-1/1,1/1} & \DriftAndRegion{-1/-1,1/0,-1/1,1/1}{(1.1,-.25)}{(-.19,1.0)}{(-.25,-1.1)}{(-1.0,-.19)}
&$1/2$
&2.6074
&$3/2$
&2.5211
&\href{https://oeis.org/A151289}{\texttt{A151289}} \\
\hline

20 &
\StepModel{-1/1,-1/-1,0/1,1/1,1/0} & \DriftAndRegion{-1/1,-1/-1,0/1,1/1,1/0}{(1.2,-.29)}{(-.21,.99)}{(-.29,-1.2)}{(-.99,-.21)} 
& $1/2$
& 2.4829
&$3/2$
&2.4829
&\href{https://oeis.org/A151317}{\texttt{A151317}} \\
\hline

25 &
\StepModel{-1/0,0/1,-1/-1,1/1,1/0} & \DriftAndRegion{-1/0,0/1,-1/-1,1/1,1/0}{(1.2,-.44)}{(-.25,.95)}{(-.44,-1.2)}{(-.95,-.25)} 
&$1/2$
&2.3976
&$3/2$
&2.3976
&\href{https://oeis.org/A151316}{\texttt{A151316}} \\
\hline

44 &
\StepModel{-1/-1,0/-1,1/0,1/1,0/1,-1/1} & \DriftAndRegion{-1/-1,0/-1,1/0,1/1,0/1,-1/1}{(.990,-.150)}{(-.131,.991)}{(-.150,-.990)}{(-.991,-.131)} 
&$1/2$
&2.6798
&$3/2$
&2.6798
&\href{https://oeis.org/A151324}{\texttt{A151324}} \\
\hline

\caption{The 5 interior-drift models with ${-\alpha_{11}}\sim 0.5$ and ${-\alpha_{10}}$ irrational.}
\label{tab:intnondfinite}
\end{longtable}

\vfill

\begin{longtable}{|M{0.8cm}|M{1.5cm}|M{1.5cm}|M{1.5cm}|M{1.5cm}|M{1.5cm}|M{1.5cm}|M{2cm}|} 
\hline
\(\#\) 
& Model 
& Drift
& $-{\alpha_{11}}$
& $-{\alpha_{10}}$
& $-{\alpha_{01}}$
& $-{\alpha_{00}}$
& OEIS\\
\hline

8 &
\StepModel{-1/0,0/1,1/0,1/-1} & \DriftAndRegion{-1/0,0/1,1/0,1/-1}{(.982,.189)}{(.231,.971)}{(.189,-.982)}{(-.971,.231)}
&1/2
&3/2
&3/2
&3.6379
&\href{https://oeis.org/A151283}{\texttt{A151283}} \\
\hline

19 &
\StepModel{-1/1,0/-1,1/1,1/-1} & \DriftAndRegion{-1/1,0/-1,1/1,1/-1}{(1.1,.26)}{(.31,.95)}{(.26,-1.1)}{(-.95,.31)} 
&1/2
&3/2
&3/2
&3.9214
&\href{https://oeis.org/A151273}{\texttt{A151273}} \\
\hline

31 &
\StepModel{-1/0,-1/1,0/1,1/-1,1/1} & \DriftAndRegion{-1/0,-1/1,0/1,1/-1,1/1}{(1.1,.31)}{(.34,.91)}{(.31,-1.1)}{(-.91,.34)} 
& 1/2
& 3/2
&3/2
&4.0704
&\href{https://oeis.org/A151304}{\texttt{A151304}} \\
\hline

39 &
\StepModel{-1/1,0/1,0/-1,1/0,1/-1} & \DriftAndRegion{-1/1,0/1,0/-1,1/0,1/-1}{(1.2,.44)}{(.48,.89)}{(.44,-1.2)}{(-.89,.48)}
& 1/2
& 3/2
& 3/2
& 4.5148
&\href{https://oeis.org/A151296}{\texttt{A151296}} \\
\hline

51 &
\StepModel{-1/0,-1/1,0/1,0/-1,1/1,1/-1} & \DriftAndRegion{-1/0,-1/1,0/1,0/-1,1/1,1/-1}{(.990,.150)}{(.177,.988)}{(.150,-.990)}{(-.988,.177)}

&1/2
&3/2
&3/2
&3.4714
&\href{https://oeis.org/A151313}{\texttt{A151313}} \\
\hline
\caption{The 5 non-polar exterior-drift models with ${-\alpha_{11}}\sim 0.5$}
\label{tab:exterior-bis-2}
\end{longtable}

\begin{longtable}{|M{0.8cm}|M{1.5cm}|M{1.5cm}|M{1.5cm}|M{1.5cm}|M{1.5cm}|M{1.5cm}|M{2cm}|} 
\hline
\(\#\) 
& Model 
& Drift
& $-{\alpha_{11}}$
& $-{\alpha_{10}}$
& $-{\alpha_{01}}$
& $-{\alpha_{00}}$
& OEIS\\
\hline

5 &

\StepModel{-1/-1,1/1,1/0,0/1} & \DriftAndRegion{-1/-1,1/1,1/0,0/1}{(1.3,-.52)}{(-.28,.95)}{(-.52,-1.3)}{(-.95,-.28)} 
& 0
& 3/2
&3/2
&2.3190
&\href{https://oeis.org/A151315}{\texttt{A151315}} \\
\hline

14 &
\StepModel{-1/0,0/1,1/1,1/-1} & 
\DriftAndRegion{-1/0,0/1,1/1,1/-1}{(1.,.12)}{(.14,1.0)}{(.12,-1.)}{(-1.0,.14)}
& 0
& 3/2
& 3/2
& 3.3876
&\href{https://oeis.org/A151293}{\texttt{A151293}} \\
\hline

21 &
\StepModel{-1/1,0/-1,0/1,1/1,1/0} & \DriftAndRegion{-1/1,0/-1,0/1,1/1,1/0}{(1.,.12)}{(.14,1.0)}{(.12,-1.)}{(-1.0,.14)}
&0
&3/2
&3/2
&3.3474
&\href{https://oeis.org/A151320}{\texttt{A151320}} \\
\hline

29 &
\StepModel{-1/1,0/1,0/-1,1/-1,1/1} & \DriftAndRegion{-1/1,0/1,0/-1,1/-1,1/1}{(1.2,.29)}{(.32,.94)}{(.29,-1.2)}{(-.94,.32)} 
&0
&3/2
&3/2
& 3.9862
&\href{https://oeis.org/A151308}{\texttt{A151308}} \\
\hline

38 &
\StepModel{-1/0,0/1,0/-1,1/0,1/1} & \DriftAndRegion{-1/0,0/1,0/-1,1/0,1/1}{(1.,-.10)}{(-.91e-1,1.0)}{(-.10,-1.)}{(-1.0,-.91e-1)} 
& 0
& 3/2
&3/2
&2.7574
&\href{https://oeis.org/A151319}{\texttt{A151319}} \\
\hline

45 &
\StepModel{-1/0,-1/1,0/1,1/1,1/0,1/-1} & \DriftAndRegion{-1/0,-1/1,0/1,1/1,1/0,1/-1}{(.965,.265)}{(.354,.934)}{(.265,-.965)}{(-.934,.354)}
& 0
& 3/2
&3/2
&4.0361
&\href{https://oeis.org/A151327}{\texttt{A151327}} \\
\hline

46 &
\StepModel{-1/-1,-1/1,0/1,1/1,1/0,1/-1} & \DriftAndRegion{-1/-1,-1/1,0/1,1/1,1/0,1/-1}{(.995,-.0130)}{(-.0129,1.00)}{(-.0130,-.995)}{(-1.00,-.0129)}
& 0
& 3/2
& 3/2
& 2.9599
&\href{https://oeis.org/A151325}{\texttt{A151325}} \\
\hline

56 &
\StepModel{-1/0,-1/1,0/1,0/-1,1/1,1/0,1/-1} & \DriftAndRegion{-1/0,-1/1,0/1,0/-1,1/1,1/0,1/-1}{(1.0,.15)}{(.16,1.0)}{(.15,-1.0)}{(-1.0,.16)}
& 0
& 3/2
& 3/2
& 3.497
&\href{https://oeis.org/A151330}{\texttt{A151330}} \\
\hline
\caption{The eight interior-drift models with ${-\alpha_{11}}\sim 0$}
\label{tab:exterior-tre}
\end{longtable} 

\end{document}